\documentclass[11pt]{article}
\usepackage{amssymb}
\usepackage{amsmath}
\usepackage{version}
\numberwithin{equation}{section}
\newtheorem{theorem}{Theorem}[section]
\newtheorem{lemma}[theorem]{Lemma}
\newtheorem{corollary}[theorem]{Corollary}

\newtheorem{definition}[theorem]{Definition}
\newtheorem{proposition}[theorem]{Proposition}

\newtheorem{note}[theorem]{Note}
\newcommand\beq{\begin{equation}}
\newcommand\ds{\displaystyle}
\newcommand\eeq{\end{equation}}
\newcommand\ii{\mathrm i}
\newcommand\al{\alpha}

\newcommand\ep{\varepsilon}
\newcommand\la{\lambda}

\newcommand{\ph}{\varphi}
\newcommand\si{\sigma}

\newcommand\D{\mathbb D}
\newcommand\T{\mathbb T}
\newcommand\C{\mathbb C}

\newcommand\R{\mathbb R}
\newcommand\e{\mathrm e}
\newcommand\Pick{\mathcal P}
\newcommand\Schur{\mathcal{S}}
\newcommand\half{{\tfrac{1}{2}}}

\newcommand\diag{\mathrm {diag~}}
\newcommand\df{\stackrel{\rm def}{=}}

\newcommand\qed{\hfill $\square$ \vspace{3mm}}
\newcommand\nn{\nonumber}
\let\phi\varphi
\DeclareMathOperator{\schur}{schur}
\DeclareMathOperator{\im}{Im}
\DeclareMathOperator{\re}{Re}

\numberwithin{equation}{section}

\title{Boundary Nevanlinna-Pick interpolation via reduction and augmentation}
\author{Jim Agler
\thanks{Partially supported by National Science Foundation Grant
DMS 0801259}
\and
N. J. Young
\thanks{Supported by EPSRC Grant EP/G000018/1} \thanks {MSC Class: 30E05}}
\date{\empty}
\begin{document}
\maketitle
\bibliographystyle{plain}

\begin{abstract}
We give an elementary proof of Sarason's solvability criterion for the Nevanlinna-Pick problem with boundary interpolation nodes and boundary target values.  We also give a concrete parametrization of all solutions of such a problem.  The proofs are based on a reduction method due to Julia and Nevanlinna. Reduction of functions corresponds to Schur complementation of the corresponding Pick matrices.
\end{abstract}
\section{Introduction} \label{intro}
In this paper we present an elementary solution of a version of the Nevanlinna-Pick interpolation problem with boundary interpolation nodes and boundary target values.  Such problems have been studied by many authors, beginning with R. Nevanlinna \cite{Nev1} in 1922, and they continue to attract interest, not only as natural questions in function theory,  but also because they have applications to engineering, particularly electrical networks and control theory (see \cite[page 812]{ball} for some references).  The problem is to construct functions in either the Schur class $\Schur$ or the Pick class $\Pick$ (defined in Section \ref{2theorems}) subject to finitely many interpolation conditions on the values of the function and its first derivative at points on the {\em boundary} of its domain of analyticity.  

 Nevanlinna, in his original paper \cite{Nev1} on  the problem, gave a simple recursive technique based on a theorem of G. Julia to determine whether a boundary interpolation problem is solvable, and if so, to describe all solutions.  The method is analogous to ``Schur reduction" for the standard Nevanlinna-Pick problem.  Nevanlinna did not, however, give a criterion for the existence of a solution.  In 1998 D. Sarason, after referring to subsequent papers by M. G. Krein, J. A. Ball, J. W. Helton, A. Kheifets and several others, wrote  \cite{S} ``It is unclear, at least to this author, to what extent the preceding treatments can be extended so as to clarify Problems $\partial NP\mathcal{S}$ and $\partial NP\mathcal{C}$, for example, and to find a criterion for the solution set to be non-empty."  He went on to give a succinct criterion, after converting the interpolation problems in question to moment problems on the circle.  He also gave a description of the solution set in terms of measures on the unit circle that solve a form of truncated trigonometric moment problem.

The elementary reduction technique of Julia and Nevanlinna has the virtue that it assumes only a good first course in complex analysis.  It does not require any operator theory or even the notion of Hilbert space, and so may be particularly suitable for working engineers.  In this paper we show that the technique can indeed be used to prove Sarason's solvability criterion, and furthermore to parametrize all solutions of Problems $\partial NP\Pick$ and $\partial NP\Schur$ (Sarason's terminology, explained below) in a straightforward and concrete way.  Nevanlinna's own parametrization \cite[Satz I]{Nev1} is marred by an oversight, which we correct in Theorem \ref{parametrize}.   We also remove an unnecessary restriction (distinctness of the interpolation values) and we incorporate the condition for indeterminateness.  It transpires that the parametrization problem is more delicate and even more interesting than Nevanlinna realized. 

There is a substantial literature on boundary interpolation.  Besides the papers discussed above we mention particularly papers of J. A. Ball and J. W. Helton \cite{bh}, D. R. Georgijevi\'c \cite{Geo} and V. Bolotnikov and A. Kheifets \cite{Bolot}, and the books of J. A. Ball, I. C. Gohberg and L. Rodman \cite{bgr}, and of V. Bolotnikov and H. Dym \cite{BD}.
These authors variously make use of Krein spaces, moment theory, measure theory, reproducing kernel theory, realization theory and de Branges space theory.  They obtain far-reaching results, including generalizations to matrix-valued functions and to functions allowed to have a limited number of poles in a disc or half plane; see \cite {bgr, Bolot} both for results and for many references to the literature.  We shall comment briefly throughout on the points of contact of their work with this paper.  An elementary approach also has its virtues, and provides different insights into the subtleties of boundary interpolation.

The paper is organised as follows.  In Section \ref{2theorems} we state the problem and  two theorems:  Sarason's solvability criterion and the corrected version of Nevanlinna's parametrization of all solutions.  We also briefly describe some other authors' parametrizations.
In Section \ref{reduction} we present an extension of the Julia-Nevanlinna reduction method and the inverse process of augmentation.  We prove the key fact that reduction and augmentation preserve the Pick class.  In Section \ref{example} we give a simple example which both illustrates Nevanlinna's procedure and provides a counterexample to \cite[Satz I]{Nev1}.   In Section \ref{schurcomp} we prove that Julia reduction of functions corresponds to Schur complementation of Pick matrices (Theorem \ref{compM}).  In Section \ref{criter} we give the promised elementary proof of Sarason's result, Theorem \ref{sarason}.  
In Section \ref{param2} we prove the parametrization theorem, Theorem \ref{parametrize}, and show that it can be interpreted as a linear fractional parametrization  or as a continued fraction expansion (Corollaries \ref{linfrac} and \ref{continued}).  In Section \ref{schur} we sketch the corresponding results for boundary interpolation in the Schur class.

We are grateful to Joseph Ball and Vladimir Bolotnikov for drawing our attention to several relevant papers.

\section{The interpolation problem and two theorems} \label{2theorems}
We denote by $\Pi$ the upper half plane $ \{ z\in\C: \im z > 0\}$, and by $\D$ the open unit disc $\{z:|z| < 1\}$.
The {\em Pick class} $\Pick$  is the set of functions  analytic and with non-negative imaginary part in $\Pi$.  The {\em Schur class} $\Schur$ is the set of functions  analytic and bounded by 1 in $\D$.  We study interpolation problems for functions in $\Pick$ and $\Schur$.   Let us state the version for $\Pick$.\\
{\noindent \bf Problem $\partial NP\Pick$:}  {\em Given distinct points $x_1, \dots, x_n \in \R$ and numbers $w_1,\dots,w_n \in \R$ and $v_1,\dots, v_n\geq 0$, determine whether there is a function $f\in\Pick$ satisfying
\beq \label{interp}
f(x_j)= w_j, \quad f'(x_j) = v_j \qquad \mbox{ for } j=1, \dots, n,
\eeq
and if there is,  describe the class of all such functions $f$.

The {\em Pick matrix} for this problem is the matrix $M= [m_{ij}]_{i,j=1}^n$ where
\[
m_{ij} = \left\{ \begin{array}{ccc} v_i & \mbox{if} & i=j \\ \ds \frac{w_i-w_j}{x_i-x_j} & \mbox{if} & i\neq j. \end{array} \right.
\]
}
For simplicity we shall understand the conditions (\ref{interp}), as Nevanlinna did, in the straightforward sense that $f$ is analytic at each interpolation node, though Sarason interpreted them in a weaker sense.   The simplification is justified by the fact that Problem $\partial NP\Pick$ has a solution in the weak sense (values in conditions (\ref{interp}) are non-tangential limits) if and only if it has a solution in the strongest possible sense (the function is rational and analytic at the interpolation nodes) -- see Note \ref{weaker} and Proposition \ref{allequiv}.  

Sarason \cite{S} studied the equivalent problems for $\Schur$ and the class $\mathcal{C}$ of analytic functions with non-negative real part in $\D$.    We state the problem $\partial NP\Schur$ in Section \ref{schur} below.

A first guess, in analogy with Pick's theorem for interpolation at interior points of $\D$ or $\Pi$, might be that there exists a solution of Problem $\partial NP\Pick$ if and only if the Pick matrix $M$ for the data is positive.  However, there is a simple reason why this cannot be so -- the fact that a non-constant function in $\Pick$ has a positive derivative at any point in $\R$ at which it takes a real value (Proposition \ref{posderiv} (1)).  Consider the following problem $\partial NP\Pick$:\\
{\em Find $f\in\Pick$ such that $f(x_1)=f(x_2)=0$ and $f'(x_1)=0, f'(x_2) = 1$, where $x_1\neq x_2\in\R$.}

The constraint $f'(x_1)=0$ forces $f$ to be constant, and so there is no such $f$, but the Pick matrix for these data is
\[
M= \left[ \begin{array}{cc} 0 & 0 \\ 0 & 1 \end{array} \right],
\]
which is positive.  One might argue that the condition $f'(x_1)=0$ is unnatural, and hope that positivity of $M$ {\em would} suffice in the case of strictly positive $v_j$.  However, it is not so:  the above example can be modified to produce a counter-example with three interpolation nodes and all $v_j>0$.

It is an intriguing feature of boundary interpolation problems that positivity of the Pick matrix is necessary but not sufficient for solvability.  As far as we know this was first pointed out by B. Abrahamse \cite[page 812]{ball}.  

Here is Sarason's solvability criterion (but for the replacement of $\mathcal{C}$ by $\Pick$).   We say that $f$ is a {\em solution} of Problem $\partial NP\Pick$ if $f\in\Pick$, $f$ is analytic at $x_1,\dots,x_n$ and $f$ satisfies the conditions (\ref{interp}).  We say that Problem $\partial NP\Pick$ is {\em solvable} if it has at least one solution, is {\em determinate} if it has exactly one solution, and is {\em indeterminate} if it has at least two solutions.  The problem has a convex solution set, and so has zero, one or infinitely many solutions.  A positive semidefinite matrix is said to be {\em minimally positive} if it majorises no non-zero positive semidefinite diagonal matrix. 
\begin{theorem} \label{sarason}
A problem $\partial NP\Pick$ 
is solvable if and only if the corresponding Pick matrix is positive definite or minimally positive.  

The solution set of the problem, if non-empty, contains a real rational function of degree no greater than the rank of the Pick matrix.

The problem is determinate if and only if the Pick matrix is minimally positive.
\end{theorem}

Our proof of Sarason's Theorem will be by means of an induction, which also leads to the following concrete parametrization of solutions.
\begin{theorem} \label{parametrize}
Suppose that Problem $\partial NP\Pick$ has a positive definite Pick matrix.  
There exist
\[
s_1,\dots,s_n\in\R\cup\{\infty\}, \qquad t_1,\dots,t_n > 0,  \qquad  y_1,\dots, y_n \in\R\cup\{\infty\}
\]
such that the general solution of Problem $\partial NP\Pick$ is $f=f_1$ where the functions $f_{n+1},\dots,f_1$ are given recursively by:
\begin{itemize}
\item[\rm(1)] $f_{n+1}$ is any function in $\Pick$ that is meromorphic at $x_1,\dots, x_n$ and satisfies $f_{n+1}(x_j) \neq y_j, \, j=1,\dots, n$;
\item[\rm (2)] $f_k$ is the augmentation of $f_{k+1}$ at $x_k$ by $s_k, t_k$ for $k=n,\dots,1$.
\end{itemize}

The quantities $s_j, t_j,y_j$ for $j=1,\dots,n$ are given explicitly by
\beq \label{defsty}
s_j = w^j_j, \qquad t_j = v^j_j, \qquad y_j=y^{n+1}_j
\eeq
where $w^k_j, v^k_j$ for $1\leq k\leq j\leq n$ and $y^k_j$ for $1\leq j < k \leq n+1$ are given by equations {\rm (\ref{defwv1})} to {\rm (\ref{defy})} below.
\end{theorem}
The notion of the augmentation of a function at a point of the real axis will be defined in the next section, and the equations {\rm (\ref{defwv1})} to {\rm (\ref{defy})} that define the quantities $s_j, t_j,y_j$ are simple piecewise rational expressions.    In this way we obtain a parametrization of  all solutions of Problem $\partial NP\Pick$ as a linear fractional transform of a function $f_{n+1}$ in the Pick class, which however is not quite free: it has to be meromorphic at the interpolation nodes and it must satisfy a single inequation at each node.   The linear fractional form of the parametrization is stated in Corollary \ref{linfrac}.  It can also be stated in continued fraction form: see Corollary \ref{continued}.

 Theorem \ref{parametrize} is a corrected version of \cite[Satz I.2]{Nev1}, where  the  inequations $f_{n+1}(x_j) \neq y_j$ were omitted on account of an oversight discussed further in Note \ref{Nserror}.     

A more general (but accordingly less explicit) result is proved in \cite[Theorems 21.1.2 and Corollary 21.4.2]{bgr} with the aid of realization theory.   In this version one seeks a rational matrix-valued function $f$  with non-negative real part and with at most $k$ poles in $\Pi$ satisfying matricial versions of the interpolation conditions (\ref{interp}).  On specialization to scalar functions we obtain the following. If the Pick matrix of Problem $\partial$NP$\Pick $ is invertible and has at most $k$ negative eigenvalues then the solution set consists of all functions of the form
\[
f(z) = (\Theta_{11}(z)g(z)+\Theta_{12}(z))(\Theta_{21}(z)g(z)+\Theta_{22}(z))^{-1}
\]
where the $\Theta_{ij}$ are given by an explicit formula in terms of the data, and $g$ ranges over the scalar rational functions that are bounded by $1$ on $\Pi$ {\em and are such that $\Theta_{21}g+\Theta_{22}$ has simple poles at the interpolation points $x_1, \dots, x_n$}.  The italicised condition is equivalent to a set of inequations $g(x_j) \neq y_j$, as in Theorem \ref{parametrize}.

Another parametrization is given in \cite[Theorem 4.3]{Bolot} within the framework of the Ukrainian school's Abstract Interpolation Problem.  The authors allow finitely many derivatives to be prescribed at each interpolation point.  The proofs depend on operator theory in de Branges-Rovnyak spaces.

J. A. Ball and J. W. Helton have developed a far-reaching theory of interpolation which they call the ``Grassmanian" approach; it makes use of the geometry of Krein spaces.  It gives a unified treatment of many classical interpolation problems, including the ones studied in this paper -- see \cite[Theorem 6.3]{bh}.  Inevitably, one pays for the generality with a less concrete parametrization.

Some terminology: for a matrix $M$, ``positive" means the same as ``positive semi-definite", and is written $M \geq 0$, whereas $M>0$ means that $M$ is positive definite.  In a block matrix 
\[
M=\left[ \begin{array}{cc} A & B \\ C & D \end{array} \right],
\]
where $A$ is non-singular, the {\em Schur complement} of $A$ is defined to be $D-CA^{-1}B$ (see for example \cite{HJ85}).  By virtue of the identity
\beq\label{schid}
M=\left[ \begin{array}{cc} A & B \\ C & D \end{array} \right] = \left[\begin{array}{cc}1 & 0 \\ CA^{-1} & 1 \end{array}\right]
 \left[\begin{array}{cc}A & 0 \\ 0 & D-CA^{-1}B \end{array}\right]
 \left[\begin{array}{cc}1 & A^{-1}B \\ 0 & 1 \end{array}\right]
\eeq
it is clear that $\mathrm{rank}~(D-CA^{-1}B) = \mathrm{rank}~M - \mathrm{rank}~A$.  Furthermore, if $A>0$ then $M\geq 0$ if and only if $C=B^*$ and $D-CA^{-1}B\geq 0$.

\section{Reduction and augmentation in the Pick class} \label{reduction}
We shall need the following basic properties of the Pick class $\Pick.$
\begin{proposition} \label{posderiv}
Let $f$ be a non-constant function in the Pick class and let $x\in\R$.
\begin{enumerate}
\item[\rm(1)] If $f$ is analytic and real-valued at $x$ then $f'(x)>0$.
\item[\rm(2)] If $f$ is meromorphic and has a pole at $x$ then $f$ has a {\em simple} pole at $x$, with a negative residue.
\end{enumerate}
\end{proposition}
\begin{proof}
(1) Suppose $f(\xi+\ii\eta)=u+\ii v$ with $\xi,\eta, u, v$ real: since $v>0$ on $\Pi$ and 
$v(x) =0$ we have  $ v_\eta (x)\geq 0$ and hence, by the Cauchy-Riemann equations, $u_\xi (x)\geq 0$. Furthermore the restriction of $v$ to a neighbourhood of $x$ in $\R$ attains its minimum  at $x$, and so $v_\xi=0$ at $x$.  Hence $f'(x)=(u_\xi+\ii v_\xi)(x) \geq 0$. 

Suppose $f'(x)=0$: we can write $f(z)=f(x)+(z-x)^\nu g(z)$ in a neighbourhood of $x$, with $g$ analytic, $g(x)\neq 0$ and $\nu \geq 2$.  Pick a neighbourhood of $x$ on which $g$ is bounded away from $0$ and $\mathrm{arg}~ g$ lies in a small interval $I$ about $\mathrm{arg}~ g(x)$.  For small enough $r>0$ and $0<\theta<\pi$, if $z=x+r \mathrm{e}^{\ii \theta}$ then
\[
\mathrm{arg}~(f(z)-f(x)) = \nu\theta + \mathrm{arg}~ g(z) \in \nu\theta+I.
\]
As $\theta$ goes from $0$ to $\pi$ the interval $\nu\theta+I$ rotates at least once around the circle, and so $\im f(z) < 0$ at some point in $\Pi$, contrary to the fact that $f\in\Pick$.  Hence $f'(x)>0$.   \\
(2)  Observe that if $f$ is a non-constant function in $\Pick$ then $-1/f \in \Pick$. If $f$ has a pole of order $\nu$ at $x$ and leading term $R/z^\nu$ then  $-1/f $ has a zero of order $\nu$ at $x$ and leading term $-z^\nu/R$.  By (1) $\nu=1$ and $-1/R > 0$.  \hfill $\Box$
\end{proof}

 Nevanlinna used a recursive technique for eliminating interpolation conditions at points on the real line.  The reduction procedure (in the case of a function analytic at an interpolation point) is due to G. Julia \cite{Ju20}.
\begin{definition} \label{defreduce} \rm
 (1)  For any non-constant function $f\in\Pick$ and any $x\in\R$ such that $f$ is analytic at $x$ and $f(x)\in\R$ we define the {\em reduction of $f$ at $x$} to be the function $g$ on $\Pi$ given by the equation
\beq \label{reducef}
g(z) = -\frac{1}{f(z)-f(x)} + \frac{1}{f'(x)(z-x)}.
\eeq
For $f\in\Pick$ having a pole at $x\in \R$ with residue $R$ we define the reduction of $f$ at $x$ to be the function $g$ on $\Pi$ given by 
\beq \label{reducepole}
g(z)= f(z) -R/(z-x).
\eeq
 (2)  For any $g\in\Pick$, any $x\in\R$ such that $g$ is meromorphic at $x$ and any $a_0\in\R\cup\{\infty\}, a_1 > 0$, we define the {\em augmentation of $g$ at $x$ by $a_0, a_1$} to be the function $f$ on $\Pi$ given by
\beq \label{augmentg}
\left\{ \begin{array}{cccll} \ds \frac{1}{f(z)-a_0} &=& \ds\frac{1}{a_1(z-x)} -g(z)&\mbox{ if }& a_0\in\R  \\ f(z) &=& g(z) - \ds\frac{1}{a_1(z-x)} & \mbox{ if } & a_0 = \infty. \end{array} \right.
\eeq
\end{definition}
The former of the equations in (\ref{augmentg}) can of course be written
\beq \label{augment2}
f(z)= a_0 + \frac{1}{\frac{1}{a_1(z-x)} -g(z)}.
\eeq
\begin{note} \label{degrees} 
\noindent \rm If $f$ is a non-constant [real] rational function in $\Pick$ of degree $N$ and $f$ is analytic at $x\in\R$ then its reduction at $x$ is a [real] rational function in 
$\Pick$ of degree at most $N-1$.

If $g\in\Pick$ is [real] rational of degree $N$ then the augmentation of $g$ at any point of analyticity in $\R$ by $a_0\in\R,\, a_1> 0$ is [real] rational of degree $N+1$, and the augmentation of $g$ at any pole in $\R$ is [real] rational of degree $N$.
\end{note}

The important property of the operations of reduction and augmentation is that they preserve the Pick class, as we now show. 
\begin{theorem} \label{propfg}
Let $x \in\R$.
\begin{enumerate}
\item[\rm(1)] If a non-constant function $f\in\Pick$  either has a pole at $x$ or is analytic and real-valued at $x$ then the reduction of $f$ at $x$ also belongs to $\Pick$ and is analytic at $x$.
\item[\rm(2)] If $g\in\Pick$ is meromorphic at $x$ and $a_0\in\R,\, a_1>0$ then the augmentation $f$ of $g$ at $x$ by $a_0,\, a_1$  belongs to $\Pick$, is analytic at $x$  and satisfies $f(x)=a_0, \, f'(x) \leq a_1.$  In fact if $g$ is analytic at $x$ then $f'(x)=a_1$, while if $g$ has a pole at $x$ of residue $R$ then $f'(x)=a_1/(1-a_1R)$ and hence $f'(x) < a_1$.
\item[\rm(3)] If $g\in\Pick$ is meromorphic at $x$ and $ a_1>0$ then the augmentation $f$ of $g$ at $x$ by $\infty,\, a_1$  belongs to $\Pick$ and has a pole at $x$  with residue $-1/a_1$.
\end{enumerate}
\end{theorem}
Julia \cite{Ju20} stated and proved necessity in (1) in the case that $f$ is analytic at $x$. Nevanlinna \cite{Nev1} gave a different proof.   E. Wigner \cite{wig},  over 30 years later, stated the analogous result for the smaller class of ``$R$ functions", that is, functions  in the Pick class that are meromorphic on the entire plane and real on the real axis.  He gave a more elementary proof, which he attributed to Schiffer and Bargmann.  Some subsequent authors have followed Wigner in calling  (1) the ``Schiffer-Bargmann lemma".  One can also give a quick proof based on Nevanlinna's integral representation of functions in the Pick class (see \cite{Donoghue}).  We shall give an elementary proof based on the idea of Schiffer and Bargmann \cite[II.8a]{wig}.

The use of augmentation at a pole will be important in our treatment.  Nevanlinna seems not to have considered it; we discuss the consequence in  Note \ref{Nserror} below and in Section \ref{param2}.  The fact that, in boundary Nevanlinna-Pick interpolation, positivity of the Pick matrix is insufficient for solvability is well explained in terms of the strict inequality ``$f'(x)<a_1$" for augmentation at a pole.

\begin{proof}
\noindent (1)  Suppose (Case (i)) that $f$ has a pole at $x$ with residue $R$.  By Proposition \ref{posderiv}, $R < 0$.  The reduction $g$ of $f$, given by equation (\ref{reducepole}), is analytic at $x$; hence there is a neighbourhood $U$ of $x$ such that $\im g$ is harmonic and continuous on $\Pi \cup U$.

Let $\ep > 0$.  On $U\cap \R$ we have $\im g = \im f \geq 0$, and so by continuity there is a neighbourhood $V(x)$ of $x$ on which $\im g \geq -\ep$.  On the other hand,
 for real $t$ distinct from $x$ (including $t=\infty$) there is a neighbourhood $V(t)$ of $t$ in $\C\cup\{\infty\}$ such that for $z\in V(t)$
\[
 \im \left( -\frac{R}{z-x}\right) \geq  -\ep
\]
and hence, for $z\in V(t)\cap \Pi$, 
\[
\im g(z) = \im f(z) - \im \frac{R}{z-x} \geq -\ep.
\]
Julia gives a geometric argument to show that $g\in\Pick$, whereas Nevanlinna concludes: it follows by a well-known application of the Poisson integral that $\im g \geq 0$ on $\Pi$.
Here is an alternative, more elementary, ending to the argument.

Suppose that $g \notin \Pick$, so that there exists $z'\in\Pi$ such that 
\[
\rho \df -\im g(z') > 0.
\]
Choose $\ep = \half \rho$.
By compactness $\R\cup\{\infty\}$ is covered by finitely many of the $V(t)$, and hence there exists $\delta>0$ such that $\im g(z) \geq -\half \rho$ for all $z$ such that $0<\im z < \delta$.   Define a contour $C=C_1 \cup C_2$ in $\Pi$ as follows.  Let $r > 2|R|/\rho$.  We take $C_1$ to be the portion of the circle about $x$ of radius $r$ that lies in the half-plane $\{\im z \geq \half \delta\}$, and $C_2$ to be the segment of the line $\{ \im z= \half\delta\}$ lying between the two intersections of the line with $C_1$.  We can assume $r$  so large and $\delta$ so small that $z'$ lies inside $C$.  Then we have $-\im g \leq \half\rho$ on $C$, but $-\im g(z') = \rho$.  This contradicts the maximum principle for harmonic functions, and so $g\in\Pick$.

Now consider Case (ii): $f$ is analytic and real-valued at $x$.  Here $f'(x) > 0$ and the reduction $g$ of $f$ at $x$ is given by equation (\ref{reducef}).  The function $-1/(f-f(x))$ is in the Pick class and has a simple pole with residue $-1/f'(x)$ at $x$.  Case (i) shows that $g\in\Pick$.

 \noindent (2) When $f$ is the  augmentation of $g$ at $x$ by $a_0\in\R, a_1>0$ we have 
\[
 -\frac{1}{f(z)-a_0} = g(z) - \frac{1}{a_1(z-x)}.
\]
The right hand side is in the Pick class and so in turn are the left hand side, $f-a_0$ and $f$.
  Furthermore, from equation (\ref{augment2}),
\beq \label{formf}
f(z)= a_0+\frac{z-x}{a_1^{-1}-(z-x)g(z)}
\eeq
Since $g$ is meromorphic at $x$ and is in $\Pick$, $g$ is either analytic at $x$ or has a simple pole at $x$ with negative residue, by Proposition \ref{posderiv}.  In either case $-(z-x)g(z)$ is analytic and non-negative-valued at $x$, and hence the denominator of the right hand side of equation (\ref{formf}) is analytic and (since $a_1^{-1}>0$) strictly positive at $x$.
  Hence $f(x)=a_0$.  If $g$ is analytic at $x$ then $ f'(x)=a_1$, while if $g$ has a pole of residue $R$ at $x$ then the denominator in equation (\ref{formf}) takes the value $a_1^{-1}-R$ at $x$, and so $f'(x)=a_1/(1-a_1R)$.    Since $a_1>0$ and $R<0$ we have $1-a_1R>1$, and so $a_1/(1-a_1R) < a_1$.\\
\noindent (3)  The statement is immediate from the second equation in (\ref{augmentg}) and the fact that the function $-1/(z-x)$ is in the Pick class.
\qed
\end{proof}
\begin{note} \label{analatx}
A function $g$ obtained by reduction is not a general element of $\Pick$.
\end{note}
For if $g$ is the reduction of a Pick function meromorphic at $x$ then $g$ is analytic at $x$.  This property propagates through a recursion to yield the $n$ inequations $f_{n+1}(x_j) \neq y_j$ in Theorem \ref{parametrize}.
\begin{note} \label{inverseops}
Reduction and augmentation at a point of analyticity are inverse operations.  
\end{note}
More precisely, 
\begin{enumerate}
\item[\rm(1)] if $f\in\Pick$ is analytic at $x$ and if $g$ is the reduction of $f$ at $x$ then $f$ is the augmentation of $g$ at $x$ by $f(x),\, f'(x)$;
\item[\rm(2)] if $g\in\Pick$ is analytic at $x$ and if $f$ is the augmentation of $g$ at $x$ by $a_0,\, a_1$ then $g$ is the reduction of $f$ at $x$.
\end{enumerate}
However, when we admit augmentation at poles, it can happen that two different functions have the same augmentation at a point.  The function $z$ is the augmentation of the zero function at $0$ by $0, 1$, and is also the augmentation of the function $-1/(2z)$ at $0$ by $0,2$.

 Since reduction preserves $\Pick$ we may consider repeated reduction of a function  $f\in\Pick$.  The process will terminate only if we reach a constant function.  This can of course happen:  if $f(z)=1/(a-z)$ for some $a\in\R$ and if $x\in \R\setminus \{a\}$ then $f\in\Pick$ and the reduction of $f$ at $x$ is the constant function $a-x$.  
\begin{note}
Reduction and augmentation also apply to a wider class of functions in $\Pick$.
\end{note}
Specifically, if $f\in\Pick$ satisfies Carath\'eodory's condition at $x\in\R$, that is,
\[
\liminf_{z\to x}\frac{\im f(z)}{\im z} < \infty,
\]
then the reduction of $f$ at $x$ exists and belongs to $\Pick$.  The augmentation $f$ of $g\in\Pick$ at $x$ by $a_0\in\R, a_1>0$ lies in $\Pick$, satisfies Carath\'eodory's condition at $x$ and $f(x)=a_0, \ f'(x)\leq a_1$, where $f'(x)$ denotes the angular derivative of $f$ at $x$.  See \cite[Section 5]{amy} for a full treatment.

\section{Two examples} \label{example}
A simple two-point interpolation problem will illustrate both the power of Nevanlinna's method and an important subtlety:

{\em  Construct all functions $f\in\Pick$, analytic at $0,1$, such that 
\beq \label{condf}
f(0)=0,\quad f(1)=1, \quad f'(0)=f'(1)=2.
\eeq }
If $f$ is such a function then we may construct the reductions $f_2,f_3$ of $f, f_2$ at $0,1$ respectively.  We find
\[
f_2(z)= -\frac{1}{f(z)} + \frac{1}{2z}, \qquad f'_2(z)= \frac{f'(z)}{f(z)^2}-\frac{1}{2z^2},
\]
so that $f_2(1)=-\half, f'_2(1)= \tfrac 32$ and
\[
f_3(z) = -\frac{1}{f_2(z) +\half} + \frac{1}{\tfrac 32 (z-1)}.
\]
By Proposition \ref{propfg} $f_2, f_3 \in\Pick$.  Furthermore, $f_3$ is analytic at $1$ and, since $f_2(0) \neq \infty$,
\beq \label{f301}
f_3(0) \neq - \tfrac 23, \qquad f_3(1) \neq \infty.
\eeq

Conversely, suppose that $f_3\in\Pick$ is meromorphic at $0,1$ and satisfies the inequations (\ref{f301}).  We can augment $f_3,f_2$ at $1, 0$ by $(-\half, \tfrac 32), (0,2)$ to obtain $f_2,f$ respectively, and we obtain $f\in\Pick$ satisfying the conditions (\ref{condf}).  We can write the general solution $f$ of our problem as a continued fraction
\[
f(z) = \frac{1}{\half z + \half -} \frac{1}{ \frac{2}{3(z-1)} - f_3(z)}
\]
parametrized by an arbitrary function $f_3\in\Pick$, meromorphic at $0,1$ and satisfying the pair of inequations (\ref{f301}).

What if we choose $f_3 \in\Pick$ that does not avoid the forbidden values -- say $f_3(1) \neq \infty$ but $f_3(0)= -\tfrac 23$?
The augmentation process still produces an $f\in\Pick$, but it is not a solution of the interpolation problem, for
\[
f(0)=0, \qquad f(1)=1, \qquad f'(1)=2 \quad \mbox{ but } f'(0) < 2.
\]
The strict inequality arises because $f_2$ has a pole at $0$, and so, according to Proposition \ref{propfg}(2), the augmentation $f$ of $f_2$ at $0$ by $0,2$ satisfies $f'(0) < 2$.

This example illustrates the need for the inequations $f_{n+1}(x_j) \neq y_j$ in Theorem \ref{parametrize} and so provides a counterexample to \cite[Satz I]{Nev1}.

A second example illustrates another fact: if one performs Nevanlinna reduction on a solvable problem $\partial NP\Pick$ with distinct target values $w_j$ one can obtain a problem in which some target values coincide, and consequently, after a second reduction, a problem in which $\infty$ is a target value.

{\em Construct all functions $f\in\Pick$ such that}
\beq\label{ex2}
f(j)=f'(j)=j, \qquad j=1,2,3.
\eeq
After reduction at $1$, we seek $f_2\in\Pick$ such that
\[
f_2(2)=f_2(3)=0, \qquad f_2'(2)=1,\quad f_2'(3)=\half.
\]
Since the target values at 2 and 3 are equal, the reduction $f_3$ of $f_2$ at 2 must satisfy:\\
$f_3$ has a pole at 3 with residue $-2$.\\

\section{Reduction and Schur complementation}\label{schurcomp}
In this section we shall show that Nevanlinna's reduction process corresponds to Schur complementation of the corresponding Pick matrices.  To this end we need to broaden the scope of the problem to encompass interpolation values at $\infty$.    Nevanlinna supposed in a footnote that the interpolation values $w_j$ were distinct, but he ignored the possibility that equal target values might nevertheless arise in the course of the recursion (as they do in example (\ref{ex2}) above). 

\noindent {\bf Problem $\partial NP\Pick_\infty$: } {\em Given distinct points $x_1, \dots, x_n \in \R$ and given $w_1, \dots,w_n \in \R \cup \{ \infty \}$ and $v_1,\dots,v_n\geq 0$ such that $v_j > 0$ whenever $w_j=\infty$, determine whether there is a function $f\in\Pick$, meromorphic at $x_1,\dots,x_n$,  such that, for $j=1,\dots, n$,
\begin{enumerate}
\item[\rm (1)]  $f(x_j)=w_j \mbox{  and  } f'(x_j)=v_j \mbox{  if  }  w_j\neq \infty, $
\item [\rm (2)] $f$ has a simple pole at $x_j$ of residue  $-1/v_j$ if $w_j=\infty$.
\end{enumerate}

The {\em Pick matrix} of the data is the $n \times n$ matrix $M=[m_{ij}]$ where
\[
m_{ij} = \left\{ \begin{array}{ccl} v_i & \mbox{ if } & 1\leq i = j \leq n\\
 \ds \frac{w_i-w_j}{x_i-x_j} & \mbox{ if }  & i\neq j, \, w_i,w_j \in \R \\
  \ds \frac{1}{x_i-x_j} &\mbox{ if } & w_i = \infty, \, w_j \in\R \\
  \ds \frac{1}{x_j-x_i} &\mbox{ if } & w_i \in\R, \, w_j = \infty \\
   0 & \mbox{ if } & w_i=w_j=\infty, \, i\neq j.
\end{array} \right.
\] }

\begin{note} \label{addinfty}
Problems $\partial NP\Pick$ and $\partial NP\Pick_\infty$ are equivalent.
\end{note}
That is, the solution sets of related problems are in bijective correspondence and the corresponding Pick matrices are diagonally congruent. 

For clearly every problem $\partial NP\Pick$ is also a $\partial NP\Pick_\infty$, with the same Pick matrix.  On the other hand,  if $f \in \Pick$ is a solution of Problem $\partial NP\Pick_\infty$ then $h \df -1/(f-\al)$, for $\al \in \R \setminus \{w_1,\dots,w_n\}$, is a solution of the problem $\partial NP\Pick$  with data
\[
x_1,\dots, x_n\in\R, \qquad w_1(\al),\dots,w_n(\al) \in\R, \qquad v_1(\al),\dots,v_n(\al) \geq 0
\]
where
\beq \label{defwjal}
w_j(\al) = \left\{ \begin{array}{ccl}  -1/(w_j-\al)&\mbox{if} & w_j\in \R \\
  0 & \mbox{if} &w_j=\infty, \end{array} \right.
\eeq
and
\beq \label{defvjal}
v_j(\al) = \left\{ \begin{array}{ccl}v_j/(w_j-\al)^2  &\mbox{if} & w_j\in \R \\
  v_j & \mbox{if} &w_j=\infty. \end{array} \right.
\eeq
Conversely, if $h\in \Pick$ satisfies the problem with these data then $f \df \al-1/h$ solves the original problem $\partial NP\Pick_\infty'$. Furthermore, if the Pick matrices of the $f$ data and the $h$ data are $M$ and $ M(\al)$ respectively then
\beq\label{MMal}
M = D  M(\al) D,
\eeq
where $D = \mathrm{diag}\{d_1,\dots,d_n\}$ and 
\beq\label{defD}
d_j=\left\{\begin{array}{cl}  w_j-\al  & \mbox{ if } w_j \in \R \\
1 & \mbox{ if } w_j=\infty. \end{array}\right.  
\eeq
Thus the Pick matrices of the two problems are diagonally congruent. \hfill $\Box$

The following statement is simple to verify.
\begin{proposition}\label{c2f}
Let $g, h$ be solutions of problems $\partial NP\Pick_\infty$ with interpolation nodes $x_1,\dots,x_n$, having corresponding Pick matrices $M_g,M_h$ respectively.  If $g=c^2h+d$ for some $c,d\in\R$ with $c\neq 0$ then
\[
M_g = \Delta M_h \Delta
\]
where $\Delta = \diag\{ c^{\ep_1}, \dots, c^{\ep_n}\}$ and
\[
\ep_j = \left\{ \begin{array}{cl} 1 & \mbox{ if } g(x_j) \in\R \\
      -1 & \mbox{ if } g(x_j) = \infty.
      \end{array} \right.
\]
\end{proposition}
  We denote by $\schur A$ the Schur complement of the {\rm (1,1)} entry of a matrix $A$, as defined at the end of the Introduction.
\begin{theorem} \label{compM}
Let Problem $\partial NP\Pick_\infty$ have Pick matrix $M$, and suppose that $n>1$ and $v_1>0$.  If $f$ is a solution of Problem $\partial NP\Pick_\infty$ then the reduction $g$ of $f$ at $x_1$ is a solution of the problem $\partial NP\Pick_\infty$ with data
\[
 x_2,\dots,x_n, \qquad w_2',\dots,w_n' \in\R\cup\{\infty\}, \qquad v_2',\dots,v_n' \geq 0
\]
and Pick matrix $\tilde M$ satisfying
\beq\label{propMtil}
\Lambda \tilde M \Lambda=   \schur M,  
\eeq
where $\Lambda =\diag\{ \la_2,\dots,\la_n\}$ and 
\beq\label{defla}
\la_j = \left\{ \begin{array}{cl} w_1 -w_j & \mbox{if } w_1,w_j\in\R\mbox{ and } w_1\neq w_j \\
       1 &   \mbox{otherwise.}
                        \end{array} \right.
\eeq
The data $w_j',v_j'$ are given explicitly by
\beq\label{defwj'}
w_j'= \left\{ \begin{array}{cll}  \ds   -\frac{1}{w_j -w_1} +  \frac{1}{v_1(x_j-x_1)} & \mbox{ if } &  w_j\neq w_1  \mbox{ in } \R \\
 \ds \frac{1}{v_1(x_j-x_1)} & \mbox{ if } &  w_1\in\R, w_j=\infty  \\
\ds  w_j+ \frac{1}{v_1(x_j-x_1)} & \mbox{ if } & w_1=\infty, w_j\in\R \\
  \infty & \mbox{ if } & w_j=w_1,
 \end{array} \right.
\eeq
\beq\label{defvj'}
v_j' = \left\{ \begin{array}{cll}\ds  \frac{v_j}{(w_j-w_1)^2} - \frac{1}{v_1(x_j-x_1)^2} & \mbox{ if } & w_j\neq w_1\mbox{ in }\R \\
 v_j-\ds \frac{1}{v_1(x_j-x_1)^2} & \mbox{ if } & w_1\in\R, w_j=\infty\mbox{ or } w_1=\infty, w_j\in\R \\
  v_j & \mbox{ if } &  w_j = w_1.
 \end{array} \right.
\eeq
\end{theorem}
\begin{proof}
The appropriate values of $w_j', v_j'$ in equations (\ref{defwj'}) and (\ref{defvj'}) are easily calculated from the relation
\[
g(z) = \left\{ \begin{array}{cll}\ds  -\frac{1}{f(z)-w_1} + \frac{1}{v_1(z-x_1)} & \mbox{ if } & w_1\in\R \\
      f(z) + \ds \frac{1}{v_1(z-x_1)} & \mbox{ if } & w_1 = \infty.  \end{array} \right.
\]
We first prove that $\Lambda \tilde M \Lambda=   \schur M$ under the assumption that
 all $w_j\in\R$.  Let
\[
M=[m_{ij}]_{i,j=1}^n, \quad \tilde M=[\tilde m_{ij}]_{i,j=2}^n, \quad\schur M = [m^\flat_{ij}]_{i,j=2}^n. 
\]
For $j=2,\dots,n$,
\[
\tilde m_{jj} = v_j' =\left\{\begin{array}{cl}\ds\frac{v_j}{(w_j-w_1)^2} - \frac{1}{v_1(x_j-x_1)^2} & \mbox{ if } w_j\neq w_1 \\
    v_j & \mbox{ if }w_j = w_1. \end{array}\right.
\]
If $w_j=w_1$ then $m_{1j}=0$ and $\la_j=1,$  while otherwise $\la_j=w_1-w_j$, and so
\[
\la_j^2 \tilde m_{jj} = m_{jj} - \ds\frac{m_{1j}m_{j1}}{m_{11}} = m_{jj}^\flat.
\]
Thus $\Lambda \tilde M \Lambda$ and  $ \schur M$ have the same diagonal entries.

For $2\leq i\neq j \leq n, \, w_i\neq w_1,\, w_j\neq w_1$,
\begin{eqnarray*}
\tilde m_{ij} &=&\ds \frac{w_i'-w_j'}{x_i-x_j} = \frac{1}{x_i-x_j}\left(-\frac{1}{w_i-w_1}+\frac{1}{v_1(x_i-x_1)} + \frac{1}{w_j-w_1} -\frac{1}{v_1(x_j-x_1)}\right)\\
   &=& \ds\frac{1}{(w_1-w_i)(w_1-w_j)} \left( m_{ij} - \frac{m_{1i}m_{j1}}{m_{11}}\right)  = \frac{m_{ij}^\flat}{\la_i \la_j}.
\end{eqnarray*}
If $w_i=w_1, w_j\neq w_1$ then $m_{i1}=0, \la_i=1$ and 
\[
\tilde m_{ij} = \ds\frac{1}{x_i-x_j}= \frac{m_{ij}}{w_i-w_j} =
   \frac{1}{w_1-w_j}\left( m_{ij} - \frac {m_{i1}m_{1j}}{m_{11}}\right) = \frac{m_{ij}^\flat}{\la_i\la_j}.
\]
A similar calculation applies if $w_i\neq w_1, w_j =w_1$, while if $w_i=w_j=w_1$ then
\[
\tilde m_{ij} = 0 = m_{ij} = \ds\frac{m_{ij}^\flat}{\la_i\la_j}.
\]
Hence $\Lambda \tilde M \Lambda = \schur M$ when all $w_j$ are finite.

Now consider the general case of Problem $\partial NP\Pick_\infty$, with some $w_j$ possibly infinite.  We reduce to the previous case by the transformation used in Note \ref{addinfty}.  Choose $\al\in\R\setminus\{w_1,\dots,w_n\}$ and consider the problem $\partial NP\Pick$ solved by $h\df -1/(f-\al)$ at $x_1,\dots,x_n$: this has data $x_j, w_j(\al)\in\R, v_j(\al)$ given by equations (\ref{defwjal}), (\ref{defvjal}), and has Pick matrix $M(\al) = D^{-1}MD^{-1}$ with $D$ as in (\ref{defD}).  Let $h_2$ be the reduction of $h$ at $x_1$.  If $\tilde M(\al)$ is the Pick matrix of $h_2$ at $x_2,\dots,x_n$, then by the previous case,
\beq\label{Mtildal}
\tilde M(\al) = \Lambda_\al^{-1}(\schur M(\al)) \Lambda_\al^{-1}
\eeq
where $\Lambda_\al = \diag \{ \la_2(\al),\dots,\la_n(\al) \}$ and
\begin{eqnarray}\label{deflaal}
\la_j(\al) &=& \left\{\begin{array}{cl} w_1(\al)-w_j(\al) & \mbox{ if }  w_1(\al)\neq w_j(\al) \\
    1 & \mbox{ otherwise} \end{array} \right.  \nn \\
   &=& \left\{\begin{array}{cl}  \ds\frac{w_1-w_j}{(w_1-\al)(w_j-\al)} & \mbox{ if } w_1,w_j \in\R, w_1 \neq w_j \\
   \ds -\frac{1}{w_1-\al} & \mbox{ if } w_1\in\R, w_j=\infty\\
   \ds \frac{1}{w_j-\al} & \mbox{ if }w_1=\infty, w_j\in\R \\
     1 & \mbox{ if } w_1=w_j.
        \end{array}\right.
\end{eqnarray}
Note that, if $D_2=\diag\{d_2,\dots,d_n\}$ then
\[
\schur M(\al) = \schur (D^{-1} M D^{-1}) = D_2^{-1}(\schur M) D_2^{-1},
\]
which, in conjunction with equation (\ref{Mtildal}),  yields
\beq\label{gotMtildal}
\tilde M(\al) = \Lambda_\al^{-1}D_2^{-1}(\schur M) D_2^{-1}\Lambda_\al^{-1}.
\eeq

A simple calculation gives the relation
\beq\label{gh2}
g = \left\{\begin{array}{cl} \ds \frac{h_2}{(w_1-\al)^2} + \frac{1}{w_1-\al} & \mbox{ if } w_1\in\R \\
   h_2+\al & \mbox{ if } w_1 = \infty. \end{array} \right.
\eeq
Thus, in the case $w_1=\infty$, 
\[
\tilde M= \tilde M(\al) = (D_2\Lambda_\al)^{-1} (\schur M) (D_2\Lambda_\al)^{-1}.
\]
On comparing equations (\ref{defD}), (\ref{deflaal}) and (\ref{defla}) we find that $D_2\Lambda_\al=\Lambda$.  Hence $\Lambda \tilde M \Lambda = \schur M$  when $w_1=\infty$.

In the case that $ w_1 \in\R$ the functions $g, h_2$ have poles at those $x_j, j=2,\dots,n$, such that $w_j=w_1$.
By equation (\ref{gh2}) and Proposition \ref{c2f},
\[
\tilde M = \Delta \tilde M(\al) \Delta
\]
where $\Delta = \diag\{(w_1-\al)^{\ep_2},\dots,(w_1-\al)^{\ep_n}\}$
and
\[
\ep_j = \left\{\begin{array}{cl} -1 & \mbox{ if } w_j \neq w_1 \\
     1 & \mbox{ if }  w_j=w_1.
     \end{array} \right.
\]
On combining this equation with (\ref{gotMtildal}) we obtain
\[
\tilde M =\Delta \Lambda_\al^{-1} D_2^{-1} (\schur M)D_2^{-1} \Lambda_\al^{-1}\Delta.
\]
By inspection,
\[
\Delta \Lambda_\al^{-1} D_2^{-1} = \Lambda^{-1},
\]
and hence $\Lambda \tilde M \Lambda = \schur M$, as required. \hfill $\Box$
\end{proof}
\begin{note} \label{detMtilde}
For $M,\tilde M$ as in Theorem {\rm \ref{compM}},
\[
\det \tilde M = \left\{\begin{array}{cl} v_1^{-1}\det M & \mbox{ if } w_1 =\infty \\
    v_1^{-1} (\det M)\prod_{w_j \in\R \setminus \{ w_1\}}(w_1-w_j)^{-2} &\mbox{ if } w_1\in\R.
   \end{array} \right.
\]
\end{note}
For
\[
\det M = v_1\det (\schur M) = v_1\det (\Lambda \tilde M \Lambda)
     = v_1(\det \Lambda)^2\det \tilde M = v_1\prod_{j=2}^n \la_j^2 (\det \tilde M).
 \]

A  result related to Theorem \ref{compM} was proved by I. V. Kovalishina \cite[Section 3]{ko89}.  She studied the analogous interpolation problem for several derivatives at a single point on the real axis.  She gave a different formula for reduction, based on a ``fundamental matrix inequality", and showed that the Pick matrix of the reduced problem is congruent to the Schur complement of the (1,1) entry of the initial Pick matrix.  Overall though there is not a great deal in common between our papers.

\section{Sarason's solvability criterion} \label{criter}
According to the Introduction of \cite{S}, the most natural analogue of the standard Nevanlinna-Pick problem for the version with boundary data   is not Problem $\partial NP\Pick$ but rather the following ``relaxed" version of it:

\noindent {\bf Problem $\partial NP\Pick'$:} {\em As Problem $\partial NP\Pick$, except that the condition $f'(x_j)=v_j$ is replaced by $f'(x_j) \leq v_j$.} 

We shall also need the relaxed version of the broadened problem, with infinities:

\noindent {\bf Problem $\partial NP\Pick'_\infty$:} {\em As Problem $\partial NP\Pick_\infty$ except that the condition {\rm (1)} $f'(x_j)=v_j$ is replaced by {\rm (1$'$)}$f'(x_j)\leq v_j$ if $w_j\in\R$, and condition {\rm (2)} is replaced by\\
{\rm ($2'$)}  $f$ has a simple pole at $x_j$ of residue at most $-1/v_j$ if $w_j=\infty$. }
\begin{note}\label{equiv'}
Problems $\partial NP\Pick'$ and $\partial NP\Pick'_\infty$ are equivalent.
\end{note}
For exactly as in Note \ref{addinfty}, the transformation $h=-1/(f-\al)$, for suitable $\al\in\R$, converts Problem $\partial NP\Pick'_\infty$ into $\partial NP\Pick'$, the two problems having diagonally congruent Pick matrices.

Sarason proves his solvability result by first obtaining one for the relaxed problem, a strategy that we shall follow.    The following inductive proof is reminiscent of D. Marshall's elementary proof of Pick's theorem via Schur reduction \cite{Mar}.

\begin{theorem} \label{probprime}
A problem $\partial NP\Pick'$ or  $\partial NP\Pick_\infty'$ is solvable if and only if the corresponding Pick matrix is positive.

The solution set of the problem, if non-empty, contains a real rational function of degree no greater than the rank of the Pick matrix.

The problem is determinate if and only if the Pick matrix is positive and singular.
\end{theorem}
\begin{proof}
In view of Note \ref{equiv'} it suffices to consider Problem $\partial NP\Pick'$.  Let $M$ be its Pick matrix.  In the case that some $v_i=0$ the problem has a solution if and only if the constant function $w_i$ is a solution,
and it is easily seen that all statements of the theorem hold.
  We may therefore assume that all $v_i$ are non-zero.

In the case $n=1$, $M=[v_1]>0$ and the problem has infinitely many solutions
\[
f(z) =  w_1+c(z-x_1) 
\]
for $0 \leq c\leq v_1$, all real rational of degree at most 1.  Thus the theorem holds when $n=1$.

Now consider $n>1$, and suppose the theorem valid for $n-1$.   Again by Note \ref{equiv'}, the version for $\partial NP\Pick_\infty'$ will also hold.

Suppose that Problem $\partial NP\Pick'$ has a solution $f$. The case in which the problem is solved by a constant function is easy, so we can assume that $f$ is non-constant, whereupon all $f'(x_j)>0$.   Consider the problem $\partial NP\Pick'$ obtained when $v_j$ is replaced by $u_j=f'(x_j)$ for $j=1,\dots,n$.
Then $u_j\leq v_j$ and the Pick matrix $M'$ of the modified problem satisfies
\beq\label{MM'}
M=M' + \diag\{v_1-u_1,\dots,v_n-u_n\} \geq M'.
\eeq
By Theorem \ref{compM}, the reduction $g$ of $f$ at $x_1$ is the solution of a problem $\partial NP\Pick'_\infty$ on $n-1$ nodes having Pick matrix $\tilde M$ satisfying
\[
\tilde M = \Lambda^{-1} (\schur M')\Lambda^{-1}
\]
for some nonsingular real diagonal matrix $\Lambda$.  By the inductive hypothesis $\tilde M \geq 0$, and hence $\schur M' \geq 0$.  It follows from the identity (\ref{schid}) that $M'\geq 0$, and hence, by the inequality (\ref{MM'}), that $M\geq0$.
 Thus the existence of a solution of $\partial NP\Pick'_\infty$ implies the positivity of the corresponding Pick matrix.

Conversely, suppose that $M\geq 0$, and let $r=\mathrm{rank}~ M$.   Then $\schur M$ is positive and of rank $r-1$.  By Theorem \ref{compM} the problem $\partial NP\Pick'_\infty$ with  data 
\[
x_2,\dots,x_n,  \quad w_2', \dots,w_n'\in\R\cup\{\infty\}, \quad v_2',\dots,v_n'>0
\]
given by equations (\ref{defwj'}),(\ref{defvj'}) has Pick matrix $\tilde M$ congruent to $\schur M$.  Thus $\tilde M \geq 0$ and $\mathrm{rank~} \tilde M = r-1$.
 By the inductive hypothesis, the reduced problem $\partial NP\Pick'_\infty$ has a solution $g\in\Pick$ that is real rational of degree at most $r-1$.   Let $f$ be the augmentation of $g$ at $x_1$ by $w_1, v_1$ (it can be that $x_1$ is a pole of $g$ -- see Note \ref{augpole} below).  By Theorem \ref{propfg} $f\in \Pick$ and $f(x_1)=w_1, \, f'(x_1) \leq v_1$, and by Note \ref{degrees}, $f$ is real rational of degree at most $r$.  We have
\[
\frac{1}{f(z)- w_1}= \frac{1}{v_1(z-x_1)} -g(z),
\]
and hence, from condition (1), if  $w_j\neq w_1$ then
\[
\frac{1}{f(x_j)-w_1} = \frac{1}{v_1(x_j-x_1)} - w_j'= \frac{1}{w_j-w_1},
\]
so that $f(x_j)=w_j$, and 
\[
\frac{f'(x_j)}{(w_j-w_1)^2}=\frac{1}{v_1(x_j-x_1)^2} + g'(x_j) \leq \frac{1}{v_1(x_j-x_1)^2}+ v_j' = \frac{v_j}{(w_j-w_1)^2}
\]
so that $f'(x_j)\leq v_j$.  If, however, $j\geq 2, w_j=w_1$ and $g$ has a simple pole of residue $R_j \leq -1/v_j$ at $x_j$ then we have
\[
\frac{1}{f(z)-w_j} = -\frac{R_j}{z-x_j} + O(1)
\]
for $z$ in a punctured neighbourhood of $x_j$, and hence $f(x_j)=w_j$ and
\[
f'(x_j) = -\frac{1}{R_j} \leq v_j.
\]
Thus $f$ is a solution of the original problem $\partial NP\Pick'$.  Hence $M\geq 0$ implies solvability of $\partial NP\Pick'$ by a real rational function of degree at most $\mathrm{rank}~ M$.  

It remains to prove the statements about determinacy.  Suppose $M$ is singular and positive.  If all the $w_j$ are equal then the off-diagonal entries of $M$ are all zero, and since $M$ is singular some $v_i$ is zero, and hence Problem $\partial NP\Pick'$ has a unique, constant, solution. We may therefore suppose that $w_1\neq w_2$.  Let $f_1, \, f_2$ be solutions of Problem $\partial NP\Pick'$.  Since $M$ is singular, so is $\tilde M$, and hence by the inductive hypothesis the problem $\partial NP\Pick'_\infty$ obtained by reduction at $x_1$ has a unique solution $g$.  Thus, for $z\in\Pi$,
\[
-\frac{1}{f_1(z)-w_1}+\frac{1}{f_1'(x_1)(z-x_1)} = g(z) = -\frac{1}{f_2(z)-w_1}+\frac{1}{f_2'(x_1)(z-x_1)}.
\]
Since $w_2\neq w_1$ we may substitute $z=x_2$ to obtain
\[
-\frac{1}{w_2-w_1}+\frac{1}{f_1'(x_1)(x_2-x_1)} = -\frac{1}{w_2-w_1}+\frac{1}{f_2'(x_1)(x_2-x_1)}.
\]
It follows in turn that $f'_1(x_1)= f'_2(x_1)$, that $-1/(f_1-w_1) = -1/(f_2-w_1)$ and that 
$f_1=f_2$.  Thus if $M$ is singular and positive then Problem $\partial NP\Pick'$ has a unique solution.  

Conversely, if $M>0$ then also $\tilde M > 0$ and so, by the inductive hypothesis the reduced problem has two distinct solutions $g_1,g_2 \in\Pick$.  The augmentations of $g_1,g_2$ at $x_1$ by $w_1, v_1$ are distinct solutions of Problem $\partial NP\Pick'$.  \hfill $\Box$
\end{proof}
\begin{note} \label{augpole}
Augmentation at poles and strict inequality in Problem $\partial NP\Pick'$.
\end{note}
In the proof of sufficiency above we construct a function $f\in\Pick$ by augmenting a rational function $g$ by $w_1,v_1$ at the interpolation node $x_1$.  In the event that $x_1$ is a pole of $g$, according to Proposition \ref{posderiv}, we have $f'(x_1) < v_1$, with strict inequality.  We saw in Section \ref{example} that one can choose  $g$ to have a pole at $x_1$.  It can even be that $g$ is unique and has a pole at $x_1$.  Consider the problem $\partial NP\Pick'$ with data
\[
\begin{array}{cccc}  j &   1  & 2  & 3\\
                 x_j & -1  & 0  &  1\\
                 w_j & -1  & 0  &  1\\
                 v_j&  1  &  2 &  1.
\end{array}
\]
Here
\[
M= \left[\begin{array}{ccc} 1& 1 & 1 \\ 1 & 2 & 1 \\ 1 & 1 & 1 \end{array} \right] \geq 0,
\qquad \mathrm{rank}~ M=2.
\]
If we reduce at $x_2=0$ we arrive at the problem of finding $g\in\Pick$ such that 
\[
  g(\pm 1) = \mp \tfrac 12, \qquad g'(\pm 1) = \tfrac 12.
\]
The unique $g\in\Pick$ satisfying these equations is $g(z) = -1/(2z)$, which does have a pole at $x_2=0$.  On augmenting $g$ at $0$ by $w_2=0, v_2=2$ we obtain the unique solution $f(z)=z$ of the problem $\partial NP\Pick'$, and as expected,
\[
 1 = f'(x_2) < v_2 = 2.
\]
\begin{note} \label{Nserror}
A statement of Nevanlinna.
\end{note}
At what point of \cite{Nev1} did Nevanlinna err?  Consider, in the light of the foregoing example, the statement  \cite[page 4]{Nev1}: 
``Notwendig und hinreichend daf\"ur, da\ss~ eine Funktion der verlangten Art existiert, ist, da\ss~  
entweder die Gleichungen (7) bestehen, in welchem Fall die Funktion eine reelle Konstante ist, oder $t_1=z_1'>0$ ist und eine Funktion $f_2(x)$ existiert, die der Bedingung (2) gen\"ugt und in den $(n-1)$ Punkten (12) die Werte (13) mit den Ableitungswerten (14) annimmt\footnote{In order that a function of the desired type exist it is necessary and sufficient that either the equations (7) hold, in which case the function is a real constant, or $t=z_1'>0$ and a function $f_2(x)$ exist which satisfies condition (2) and, at the $(n-1)$ points (12), takes the values (13) with corresponding derivative values (14).}."

In the terminology of this paper, he is asserting that Problem $\partial NP\Pick$ is solvable if and only if either it is solved by a constant function or its reduction at an interpolation node is solvable.  He evidently overlooked the possibility described in Note \ref{augpole}: in proving sufficiency one  augments $f_2$ at $x_1$, and $x_1$ may be a pole of $f_2$.  In this case one obtains a strict inequality: $f'(x_1)< z_1'$.  Indeed, immediately before the above paragraph he wrote ``Ferner best\"atigt man unmittelbar, da\ss~  die Werte der Funktion und ihrer Ableitung in den Punkten (3) mit den gegebenen Werten (4) bzw. (5) \"ubereinstimmen
\footnote{Furthermore, one sees immediately that the values of the function and its derivative at the points (3) coincide with the prescribed values (4) and (5) respectively.}," a claim which fails in the case of augmentation at a pole. \hfill $\Box$

We now deduce Sarason's solvability theorem from Theorem \ref{probprime}. 

\noindent {\bf Proof of Theorem \ref{sarason}}
Necessity.  Suppose that Problem $\partial NP\Pick$ has a solution $f\in\Pick$ but that its  Pick matrix $M$ is neither positive definite nor minimally positive.  {\em A fortiori} $f$ solves Problem $\partial NP\Pick'$, and so, by Theorem \ref{probprime}, $M\geq 0$.  Since $M$ is not positive definite, $M$ is singular, and so Problem $\partial NP\Pick'$ has the {\em unique} solution $f$.  Since $M$ is not minimally positive there is some index $i$ and some positive $v_i{'} < v_i$ such that $M'\geq 0$, where $M'$ is the matrix obtained when the $(i,i)$ entry $v_i$ of $M$ is replaced by $v_i{'}$.  Again by Theorem \ref{probprime}, there exists $h\in\Pick$ such that $h(x_j)=w_j$ for each $j$, $h'(x_j)\leq v_j$ for $j\neq i$ and $h'(x_i) \leq v_i{'}<v_i$.  In view of the last relation and the fact that $f'(x_i)=v_i$ we have $h\neq f$, while clearly $h$ is a solution of Problem $\partial NP\Pick'$, as is $f$.  This contradicts the uniqueness of the solution $f$.  Hence if the problem is solvable then either $M>0$ or $M$ is minimally positive.

Sufficiency.  Suppose that $M$ is minimally positive.  By Theorem \ref{probprime} there is an $f\in\Pick$ such that $f(x_j)=w_j$ and $f'(x_j)\leq v_j$ for $j=1,\dots,n$.
If in fact $f'(x_i)<v_i$ for some index $i$ then consider the matrix $M'$ obtained when the $(i,i)$ entry $v_i$ of $M$ is replaced by $f'(x_i)$.  Since $M'$ is the Pick matrix of a problem $\partial NP\Pick'$ that is solvable (by $f$), we have $M'\geq 0$ by Theorem \ref{probprime}, and so $M$ majorises the non-zero positive diagonal matrix $\mathrm{diag}~\{0,\dots, v_i-f'(x_i), \dots,0\}$, contrary to hypothesis.  Thus $f'(x_j)=v_j$ for each $j$, that is, $f$ is a solution of Problem $\partial NP\Pick$.  Moreover, since $M$ is singular, $f$ is the unique solution of Problem $\partial NP\Pick'$, hence is real rational of degree at most $\mathrm{rank}~M$.

The above two arguments are essentially Sarason's.  For the case $M>0$ we give a new, elementary, argument based on Nevanlinna reduction.  We prove by induction a strengthening of the statement of the theorem.  Let $S(n)$ be the assertion:\\
{\em If the Pick matrix of Problem $\partial NP\Pick$ is positive definite, if $T$ is a finite subset of $\R\setminus\{x_1,\dots,x_n\}$ and $y_t\in\R\cup\{\infty\}$ for every $t\in T$ then there are at least two functions $f\in\Pick$ that solve Problem $\partial NP\Pick$,  are real rational of degree at most $n$ and satisfy $f(t)\neq y_t$ for every $t\in T$. }

The point of the condition $f(t) \neq y_t$ is to ensure that all augmentations take place at points of analyticity.  Once $S(n)$ is established we can apply it with $T$ the empty set to conclude the proof of the theorem.

$S(n)$ is true if $n=0$. Here the Pick matrix condition is vacuously true, and one may choose any two distinct real-valued constant functions $f$ with values outside the finite set $\{y_t: t \in T\}$.

Let $n\geq 1$ and suppose that $S(n-1)$ is true.  By the equivalence of Problems $\partial NP\Pick$ and $\partial NP\Pick_\infty$ the analogous statement will be true for $\partial NP\Pick_\infty$.  We apply reduction at $x_1$.  Consider the problem $\partial NP\Pick_\infty$ with data 
\[
x_2,\dots,x_n, \quad w_2',\dots,w_n', \quad v_2',\dots, v_n'>0
\]
as in Theorem \ref{compM}.
Let $T^1= T\cup\{x_1\}$, let $y^1_{x_1}=\infty $ and for $t\in T$ let
\[
y_t^1=-\frac{1}{y_t-w_1} + \frac{1}{v_1(t-x_1)},
\]
with the natural interpretation $y_t^1=\infty$ if $y_t=w_1$.
By Theorem \ref{compM}, the Pick matrix $\tilde M$ of this reduced problem is congruent to  $\schur M$, hence is positive definite.
By the inductive hypothesis there are at least two real rational functions $g\in\Pick$ of degree at most $n-1$ that solve the reduced problem $\partial NP\Pick_\infty$, are meromorphic at every point of $T^1$ and satisfy $g(t)\neq y_t^1$ for $t\in T^1$.  In particular, $g$ is analytic at $x_1$.

For either $g$, let $f$ be the augmentation of $g$ at $x_1$ by $w_1,\, v_1$:
\beq\label{gotf}
f(z) = w_1 + \frac{1}{\frac{1}{v_1(z-x_1)}-g(z)}.
\eeq
  Since $g$ is analytic at $x_1$,  Theorem \ref{propfg} and Note \ref{degrees} tell us that $f$ is analytic at $x_1$, $f(x_1)=w_1, \, f'(x_1)=v_1$ and $f$ is real rational of degree at most $n$.  For $j=2,\dots,n$ we have
\[
g(x_j)=w_j'=\left\{\begin{array}{cl} \ds -\frac{1}{w_j-w_1}+\frac{1}{v_1(x_j-x_1)} &    
 \mbox{ if } w_j\neq w_1 \\
   \infty & \mbox{ if } w_j=w_1. \end{array} \right.
\]
If $w_j\neq w_1$ then the denominator in the right hand side of equation (\ref{gotf}) is non-zero at $x_j$ and one sees that $f$ is analytic at $x_j$ and that $f(x_j)=w_j, \, f'(x_j)=v_j$.  These relations also hold if $w_j=w_1$.

We claim that $f$ is analytic on $T$ and satisfies $f(t)\neq y_t, \, t\in T$.
Indeed, since $g$ is analytic at $t$ and $g(t)\neq y^1_t$, we have 
\[
\frac{1}{v_1(z-x_1)} - g(z)
\]
is analytic and non-zero at $t$.  Hence $f$ is analytic at $t$.  If 
 $y_t\neq w_1$, then since $g(t) \neq y^1_t$ we have
\[
g(t) \neq -\frac{1}{y_t-w_1}+\frac{1}{v_1(t-x_1)}
\]
and hence
\[
\frac{1}{f(t)-w_1}=\frac{1}{v_1(t-x_1)}-g(t) \neq  \frac{1}{y_t-w_1},
\]
and therefore $f(t)\neq y_t$.  On the other hand, if it happens that $w_1=y_t$, since $g$ is analytic at $t$, it is clear that $1/(f-w_1)$ is analytic at $t$, hence $f(t)\neq w_1$.  
Thus $f(t) \neq y_t$.

In this way we construct at least two functions in $\Pick$ with the required properties.  It follows by induction that $S(n)$ holds for all non-negative integers $n$.
\hfill $\Box$

There are other solvability criteria.  Sarason gives a second one himself  \cite[Theorem 2]{S} based on an idea of Kotelyanskii, and D. R. Georgijevi\'c \cite[Theorem 1]{Geo} gave another (for a more general problem) at about the same time.  An earlier criterion may be found in \cite[Theorem XIII.1]{Donoghue}.  However,  Sarason's first criterion, as in Theorem \ref{sarason}, is the simplest and most elegant.
\begin{note}
Solvability of Problem $\partial NP\Pick_\infty$.
\end{note}
The above proof, with minimal changes, shows that Theorem \ref{sarason} remains true if 
$\partial NP\Pick$ is replaced by $\partial NP\Pick_\infty$.
\begin{note} \label{weaker}
Functions not analytic at the interpolation nodes.
\end{note}
Sarason and other authors (e.g. \cite{Bolot}) study problems like $\partial NP\Pick'$ with less stringent regularity assumptions than analyticity at the interpolation nodes.  Let us say that a function $f\in\Pick$ is a {\em weak solution} of Problem $\partial NP\Pick'$ if, for $1\leq j\leq n$, $f$ has non-tangential limit $f(x_j)=w_j$ and angular derivative $f'(x_j)$ which is no greater than $v_j$ at $x_j$.  
\begin{proposition} \label{allequiv}
The following five statements are equivalent:
\begin{enumerate}
\item[\rm(1)]  Problem $\partial NP\Pick'$ has a weak solution;
\item[\rm(2)]  Problem $\partial NP\Pick'$ has a solution;
\item[\rm(3)]  Problem $\partial NP\Pick'$ has a rational solution;
\item[\rm(4)]  Problem $\partial NP\Pick'$ has a real rational solution;
\item[\rm(5)]  The Pick matrix of the problem is positive.
\end{enumerate} 
\end{proposition}
\begin{proof}  By Theorem \ref{probprime}, (2)$\Leftrightarrow$ (3) $\Leftrightarrow$ (4) $\Leftrightarrow$ (5), and obviously (2)$\Rightarrow$(1).  It is also true that (1)$\Rightarrow$(5).  For suppose that $f$ is a weak solution of Problem $\partial NP\Pick'$. 
For any $\ep>0$ let $M(\ep)=[m_{ij}(\ep)]$ where
\[
m_{ij}(\ep) = \frac{f(x_i+\ii\ep) - \overline{f(x_j+\ii\ep)}}{x_i+\ii\ep -(x_j-\ii\ep)}
\]
for $1\leq i,j\leq n$.
By the original theorem of Pick, $M(\ep) \geq 0$.  By the Carath\'eodory-Julia Theorem (for example, \cite[VI]{S2}), as $\ep\to 0$,
\[
m_{jj}(\ep)=\frac{2\ii\im f(x_j+\ii\ep)}{2\ii\ep} \to f'(x_j). 
\]
It follows that $M(\ep)$ tends to a matrix $M'\geq 0$ which is related to the Pick matrix $M$ of Problem $\partial NP\Pick'$ by
\[
M = M' + \diag\{v_1-f'(x_1),\dots,v_n -f'(x_n)\}.
\]
Since each $f'(x_j)\leq v_j$ we have $M\geq 0$.  Thus (1)$\Rightarrow$(5). \hfill $\Box$
\end{proof}

In view of these equivalences  there is no loss (as far as existence goes) in restricting ourselves to interpolation by functions that are analytic at the interpolation nodes.

\section{Parametrization of all solutions} \label{param2}
We turn to the proof of Theorem \ref{parametrize}.  The method of proof is essentially that of Nevanlinna, but in adding the necessary considerations of infinities and augmentation at poles one cannot avoid a modicum of algebraic detail.
\begin{lemma} \label{inductf}
Let  the problem  $\partial NP \Pick_\infty$ with data
\[
x_1,\dots,x_n\in\R,\qquad 
w_1,\dots,w_n\in\R\cup\{\infty\}, \qquad v_1,\dots,v_n > 0
\]
have positive definite Pick matrix.   Define
\[
 w^k_k,\dots, w^k_n \in \R \cup \{\infty\}, \qquad v^k_k,\dots,v^k_n > 0, \quad \mbox{ for } k=1,\dots,n,
\]
recursively by the relations 
\beq \label{defwv1}
w^1_j = w_j, \qquad v^1_j = v_j, \qquad \mbox{for } j=1,\dots,n
\eeq
and, for $2\leq k+1 \leq j\leq n$,
\beq \label{defwkj}
w^{k+1}_j = \left\{ \begin{array}{ccl} 
\ds \frac{1}{v^k_k(x_j-x_k)} & \mbox{ if } & w^k_k \in \R
\mbox{ and }w^k_j=\infty\\
\ds -\frac{1}{w^k_j-w^k_k} + \frac{1}{v^k_k(x_j-x_k)} & \mbox{ if } & w^k_k,w^k_j \in \R \mbox{ and } w^k_j \neq w^k_k \\
\ds w^k_j + \frac{1}{v^k_k(x_j-x_k)} & \mbox{ if } & w^k_k =\infty\mbox{ and }w^k_j\in\R\\
\infty & \mbox{ if } &  w^k_k=w^k_j,
 \end{array}
\right.
\eeq
\beq \label {defvkj}
v^{k+1}_j = \left\{\begin{array}{ccl} 
\ds v^k_j - \frac{1}{v^k_k(x_j-x_k)^2} & \mbox{ if } & w^k_k\in\R, w^k_j=\infty \mbox{ or } w^k_k=\infty, w^k_j\in\R \\
\ds \frac{v^k_j}{(w^k_j-w^k_k)^2} - \frac{1}{v^k_k(x_j-x_k)^2} &\mbox{ if } &  w^k_k,w^k_j\in \R \mbox{ and } w^k_j \neq w^k_k \\
  v^k_j  &  \mbox{ if } & w^k_k=w^k_j. 
\end{array}\right.
\eeq
Suppose further that $y^k_1,\dots,y^k_{k-1} \in\R\cup \{\infty\}$ are defined for $k=2,\dots,n+1$ by $y^{k}_{k-1} = \infty$
and, for $1\leq j < k \leq n$,
\beq \label{defy}
y^{k+1}_j = \left\{ \begin{array}{cll}\ds -\frac{1}{y^k_j-w^k_k}+\frac{1}{v^k_k(x_j - x_k)}  &\mbox{ if } & w^k_k \in\R \\
 & & \\
  \ds y^k_j+  \frac{1}{v^k_k(x_j - x_k)}    &\mbox{ if } & w^k_k = \infty. \end{array}\right.
\eeq

For $k=1,\dots,n$,  if a function $f\in\Pick$ satisfies
\begin{enumerate}
\item[\rm(1)] $f$ is a solution of problem $\partial NP\Pick_\infty$ with data
\beq \label{NPk}
x_k,\dots,x_n \in\R, \qquad w^k_k,\dots,w^k_n \in \R\cup\{\infty\}, \qquad v^k_k, \dots, v^k_n > 0,
\eeq
\item[\rm(2)] $f$ is meromorphic at $x_1,\dots,x_{k-1}$ and $f(x_j) \neq y^k_j$ for $j=1,\dots,k-1$
\end{enumerate}
then the reduction $g$ of $f$ at $x_k$ satisfies
\begin{enumerate}
\item[\rm(1$'$)] $g$ is a solution of problem $\partial NP\Pick_\infty$ with data
\[
x_{k+1},\dots,x_n \in\R, \qquad w^{k+1}_{k+1},\dots,w^{k+1}_n \in \R\cup\{\infty\}, \qquad v^{k+1}_{k+1}, \dots, v^{k+1}_n > 0,
\]
\item[\rm(2$'$)] $g$ is meromorphic at $x_1,\dots,x_{k}$ and $g(x_j) \neq y^{k+1}_j$ for $j=1,\dots,k$.
\end{enumerate}
Conversely, if $g\in\Pick$ satisfies conditions {\rm(1$'$)} and {\rm(2$'$)} then the augmentation $f$ of $g$ at $x_k$ by $w^k_k, v^k_k$ satisfies conditions {\rm (1)} and {\rm (2)}.
\end{lemma}
The expressions on the right hand side of equation (\ref{defy}) have their natural interpretations in the event that $y^k_j=\infty$ or $y^k_j=w^k_k$.

\begin{proof}
The proof is a re-run of the proof by induction of $S(n)$ in the final part of the proof of Theorem \ref{sarason}, but with extra detail. 
If $M_k$ is the Pick matrix of the interpolation problem with data (\ref{NPk}), $k=1,\dots,n$, then for $k\leq n-1$, $M_{k+1}$ is  congruent to $\schur M_k$.  Since $M_1 > 0$ it follows that $M_k>0$ for $k=1,\dots,n$.  In particular $v^k_j > 0$ for $1\leq k\leq j\leq n$.

Fix $k$, $1\leq k\leq n$.   Suppose that $f\in\Pick$ satisfies conditions (1) and (2) and that $g$ is the reduction of $f$ at $x_k$.  By Theorem \ref{propfg}, $g\in\Pick$ and $g$ is analytic at $x_k$.  We shall show that $g$ satisfies (1$'$) and (2$'$).

Consider first the case that $w^k_k= \infty$.  Since $f$ has a pole at $x_k$ with residue $-1/v^k_k$  we have
\[
g(z) = f(z) + \frac{1}{v^k_k(z-x_k)}.
\]
Clearly $g$ is meromorphic at $x_1, \dots,x_n$.  Consider $g(x_j)$ for $j= k+1,\dots,n$.
If $w^k_j \in\R$ then $f(x_j) = w^k_j$ and so
\[
g(x_j) = w^k_j + \frac{1}{v^k_k(x_j-x_k)} = w^{k+1}_j,
\]
while if $w^k_j=\infty$ then $f$ has a simple pole at $x_j$ with residue $-1/v^k_j$, and hence the same is true for $g$.  Since $v^{k+1}_j = v^k_j$ when $w^k_k=w^k_j=\infty$ it follows that $g$ has a simple pole with residue $-1/v^{k+1}_j$ at $x_j$.  Hence $g$ satisfies (1$'$) when $w^k_k=\infty$.

Since $g$ is analytic at $x_k$ we have $g(x_k) \neq \infty = y^{k+1}_k$.  For $j=1,\dots, k-1$ we have $f(x_j) \neq y^k_j$ and so
\[
g(x_j) \neq y^k_j + \ds \frac{1}{v^k_k(x_j-x_k)} = y^{k+1}_j
\]
(note that this assertion remains valid when $y^k_j= \infty$).  Thus $g$ satisfies (1$'$) and (2$'$) when $w^k_k=\infty$.

Now consider the case that $w^k_k \in\R$.  We have
\beq \label{formg}
g(z)= \ds -\frac{1}{f(z)-w^k_k} + \frac{1}{v^k_k(z-x_k)}.
\eeq
For $j=k+1, \dots, n$ we consider the subcases (i) $w^k_j\in\R, w^k_j\neq w^k_k$, (ii) $w^k_j=w^k_k \in\R$, (iii) $w^k_j=\infty$.   In subcase (i) we have
\begin{eqnarray*}
g(x_j) &=& \ds -\frac{1}{w^k_j-w^k_k} + \frac{1}{v^k_k(x_j-x_k)} = w^{k+1}_j, \\
g'(x_j) &=& \ds \frac{v^k_j}{(w^k_j-w^k_k)^2}  - \frac{1}{v^k_k(x_j-x_k)^2}= v^{k+1}_j
\end{eqnarray*}
and so $g$ satisfies (1$'$) at $x_j$.  In subcase (ii) one sees from equation (\ref{formg}) that $g$ has a simple pole at $x_j$ with residue $-1/v^k_j$.  Since in this subcase $w^{k+1}_j=\infty, v^{k+1}_j=v^k_j$ it is again true that (1$'$) holds at $x_j$.  In subcase (iii) $f$ has a pole of residue $-1/v^k_j$ at $x_j$, and we find that
\[
g(x_j)= \ds \frac{1}{v^k_k(x_j-x_k)} = w^{k+1}_j, \quad
g'(x_j) = \ds v^k_j - \frac{1}{v^k_k(x_j-x_k)^2} = v^{k+1}_j.
\]
Hence $g$ satisfies (1$'$) when $w^k_k\in\R$.

We must show that $g$ satisfies (2$'$) when $w^k_k\in\R$.   Again $g$ is meromorphic at $x_1,\dots,x_k$ and $g(x_k) \neq \infty = y^{k+1}_k$.  For $j=1,\dots,k-1$ we have by hypothesis $f(x_j) \neq y^k_j$.  If $y^k_j \neq w^k_k$ then 
\[
g(x_j) \neq \ds -\frac{1}{y^k_j-w^k_k} + \frac{1}{v^k_k(x_j-x_k)} = y^{k+1}_j,
\]
while if $y^k_j = w^k_k$ then $g(x_j) \neq \infty =y^{k+1}_j$.
Thus $g$ satisfies (2$'$) as required.  We have shown that in all cases $g$ satisfies both (1$'$) and (2$'$).

Conversely, suppose that $g\in\Pick$ satisfies conditions (1$'$) and (2$'$), and let $f$ be the augmentation of $g$ at $x_k$ by $w^k_k, v^k_k$.  We shall prove that $f$ satisfies conditions (1) and (2).

Consider the case that $w^k_k=\infty$.  We have
\[
f(z) = g(z) - \ds \frac{1}{v^k_k(z-x_k)}.
\]
Clearly $f$ is meromorphic at $x_1,\dots,x_n$.  Since $g(x_k)\neq y^{k+1}_k=\infty$, $g$ is analytic at $x_k$.  Hence $f$ has a pole at $x_k$ with residue $-1/v^k_k$, so that $f$ satisfies (1) at $x_k$.  For $j=k+1,\dots,n$ we have
\beq \label{gxj}
g(x_j) = w^{k+1}_j.
\eeq
Either (i) $w^{k+1}_j =\infty$ or (ii) $w^{k+1}_j \in\R$.  Case (i) occurs when $w^k_k=w^k_j$.   Here $f,g$ both have poles at $x_j$ with residue $-1/v^{k+1}_j$.   Since $v^{k+1}_j=v^k_j$ when $w^k_k=w^k_j$ we deduce that $f$ has a pole at $x_j$ with residue $-1/v^k_j$.  In Case (ii)
\begin{eqnarray*}
f(x_j)&=& w^{k+1}_j - \ds \frac{1}{v^k_k(x_j-x_k)} = w^k_j \in\R, \\
f'(x_j) &=& v^{k+1}_j+\ds\frac{1}{v^k_k(x_j-x_k)^2} = v^k_j.
\end{eqnarray*}
Thus $f$ satisfies condition (1) when $w^k_k =\infty$.

For $j= 1,\dots,k-1$ we have $g(x_j) \neq y^{k+1}_j$.  Either (i) $y^{k+1}_j =\infty$ or (ii) $y^{k+1}_j \in\R$.   In Case (i) $g$ is analytic at $x_j$, and hence $f$ is also, whence $f(x_j) \neq \infty = y^k_j$.  In Case (ii)
\[
f(x_j) \neq y^{k+1}_j - \ds \frac {1}{v^k_k(x_j-x_k)} = y^k_j.
\]
Thus $f$ satisfies (1) and (2) when $w^k_k=\infty$.

The remaining possibility is that $w^k_k \in\R$.  Here 
\beq \label{fk}
\ds \frac{1}{f(z)-w^k_k} = \frac{1}{v^k_k(z-x_k)} - g(z).
\eeq
Since $g$ is analytic at $x_k$ we have $f(x_k)=w^k_k, \, f'(x_k) = v^k_k$ by Theorem \ref{propfg}.   For $j= k+1,\dots,n$ either (i) $w^{k+1}_j = \infty$ or (ii) $w^{k+1}_j \in\R$.  In Case (i) $g$ has a pole at $x_j$ with residue $-1/v^{k+1}_j$, and equation(\ref{fk}) shows that $f(x_j)= w^k_k, \, f'(x_j) = v^{k+1}_j$.
It is clear from equations (\ref{defwkj}) and (\ref{defvkj}) that $w^k_j=w^k_k$ and $v^k_j= v^{k+1}_j$, and so
\[
f(x_j)= w^k_j, \quad  f'(x_j) = v^{k}_j.
\]
In Case (ii) one of the first two alternatives in equation (\ref{defwkj}) holds, and so either 
\[
 w^k_j=\infty, \quad w^{k+1}_j = \ds \frac{1}{v^k_k(x_j-x_k)}, \quad
v^{k+1}_j=v^k-\frac{1}{v^k_k(x_j-x_k)^2}
\]
or
\[
w^k_j \in\R, w^k_j\neq w^k_k, \quad w^{k+1}_j = -\frac{1}{w^k_j-w^k_k} + \frac{1}{v^k_k(x_j-x_k)}, \quad v^{k+1}_j=\frac{v^k_k}{(w^k_j-w^k_k)^2} - \frac{1}{v^k_k(x_j-x_k)^2}.
\]
In the former case it is clear from equation (\ref{fk}) that $1/(f-w^k_k)$ vanishes  and 
has derivative
\[
\ds -\frac{1}{v^k_k(x_j-x_k)^2} - v^{k+1}_j = - v^k_j
\]
at $x_j$, which implies that $f$ has a pole at $x_j$ of residue $ -1/v^k_j$.

In the latter case ($w^k_j\in\R,\, w^k_j \neq w^k_k)$,
\[
\ds \frac{1}{f(x_j)-w^k_j }= \frac{1}{v^k_k(x_j-x_k)} - w^{k+1}_j = \frac{1}{w^k_j-w^k_k}
\]
and so $f(x_j)=w^k_j$.  On differentiating equation (\ref{fk}) and setting $z=x_j$ we have
\[
\ds -\frac{f'(x_j)}{(w^k_j-w^k_k)^2} = -\frac{1}{v^k_k(x_j-x_k)^2} - v^{k+1}_j = -\frac{v^k_j}{(w^k_j-w^k_k)^2}
\]
and so $f'(x_j)=v^k_j$.  Thus in all cases $f$ satisfies (1).

We must show that $f$ satisfies (2) when $w^k_k\in\R$.  For $j=1,\dots,k-1$ we have $g(x_j) \neq y^{k+1}_j$ and
\beq \label{dagger}
\ds \frac{1}{f(x_j)-w^k_k} = \frac{1}{v^k_k(x_j-x_k)} - g(x_j) \neq \frac{1}{v^k_k(x_j-x_k)} - y^{k+1}_j.
\eeq
If $y^{k+1}_j=\infty$ then $y^k_j=w^k_k$.  Since $1/(f-w^k_k) \neq \infty$ at $x_j$ we have $f(x_j) \neq w^k_k = y^k_k$.  On the other hand, if $y^{k+1}_j\in\R$ then equation (\ref{dagger}) becomes
\[
\ds\frac{1}{f(x_j)-w^k_k} \neq \frac{1}{y^k_j-w^k_k},
\]
from which it follows that $f(x_j) \neq y^k_j$.  This concludes the proof that $f$ satisfies conditions (1) and (2). \hfill $\Box$
\end{proof}

Note that if one starts the recursion of Theorem \ref{parametrize} with a function $f_{n+1} \in\Pick$ that is meromorphic at each $x_j$ but does {\em not} satisfy the inequations $f_{n+1}(x_j) \neq y_j$ then one obtains a function $f\in\Pick$ that is a solution of Problem $\partial NP\Pick'$, not $\partial NP\Pick$.

We can give an alternative expression for the parametrization of solutions in Theorem \ref{parametrize} in terms of a linear fractional transformation (as in \cite[\S 4]{Nev1}, where Nevanlinna considers the generic case, in which all $s_k$ are finite).  For any $2\times 2$ matrix $A=[a_{ij}]$ let us denote the corresponding linear fractional transformation by $L[A]$:
\[
L[A](w) = \frac{a_{11} w+a_{12}}{a_{21}w + a_{22}}, \qquad w\in\C\cup \{\infty\}.
\]
In this notation the relationship between $f_{k+1}$ and its augmentation $f_k$ at $x_k$ by $s_k, t_k$ can be written
\[
f_k(z) = L[A_k(z)] ( f_{k+1}(z))
\]
where, for $k=1,\dots, n,$
\beq \label{defA}
A_k(z) = \left\{ \begin{array}{lll}
  \left[\begin{array}{cc} t_k(z-x_k) & 1 \\ 0 & t_k(z-x_k)\end{array}\right] & \mbox{ if } & s_k=\infty \\
  & & \\
  \left[\begin{array}{cc} s_k t_k(z-x_k) & -t_k(z-x_k)-s_k \\  t_k(z-x_k) & -1 \end{array}\right] & \mbox{ if } & s_k\in\R. \end{array} \right .
\eeq
Note that, in either case, $\det A_k(z) = t^2_k(z-x_k)^2$.
The recursion for $f$ in Theorem \ref{parametrize} becomes:
\begin{eqnarray*}
f(z)&=& f_1(z) = L[A_1(z)](f_2(z)) = L[A_1(z)A_2(z)](f_3(z)) = \dots \\
   &=& L[A_1(z)\dots A_n(z)](f_{n+1}(z)).
\end{eqnarray*}
We therefore arrive at the following linear fractional parametrization.
\begin{corollary} \label{linfrac}
If Problem $\partial NP\Pick$ has a positive definite Pick matrix then its general solution is
\[
f(z)= \ds \frac  {a(z) h(z) + b(z)}{c(z) h(z) + d(z)}
\]
where $a, b, c, d$ are real polynomials of degree at most $n$ satisfying, for some $m>0$,
\[
(ad-bc)(z) =  m\prod_{k=1}^n  (z-x_k)^2
\]
and given by
\[
\left[\begin{array}{cc} a(z) & b(z) \\ c(z) & d(z) \end{array}\right] =
A_1(z) A_2(z) \dots A_n(z),
\]
where $A_k(z)$ is given by equations {\rm (\ref{defA})}, the quantities $s_k,t_k$ and $y_k$ are as in Theorem {\rm \ref{parametrize}} and $h$ is an arbitrary function in the Pick class that is meromorphic at $x_1,\dots,x_n$ and satisfies $h(x_j) \neq y_j, \, j=1,\dots,n$.
\end{corollary}
Yet another way to write down the parametrization of solutions in Theorem \ref{parametrize} is as a continued fraction.  In the event that $s_1,\dots,s_n$ are all finite we have, for $1\leq k\leq n$,
\[
f_k(z) = s_k + \frac{1}{\frac{1}{t_k(z-x_k)} - f_{k+1}(z)},
\]
and so
\begin{eqnarray*}
f(z) &=& f_1(z) = s_1+  \frac{1}{\frac{1}{t_1(z-x_1)} - f_{2}(z)}\\
  &=& s_1+  \frac{1}{\frac{1}{t_1(z-x_1)} -s_2- }\frac{1}{\frac{1}{t_2(z-x_2)}-f_3(z)}.
\end{eqnarray*}
Continuing in this fashion we obtain the following statement.
\begin{corollary}\label{continued}
Let Problem $\partial NP\Pick$ have a positive definite Pick matrix and let $s_k, t_k$ and $y_k$ be defined as in Theorem {\rm\ref{parametrize}} for $k=1,\dots,n$.  Suppose that $s_1,\dots,s_n \in\R$.
The general solution of the problem is
\[
f(z) = s_1+  \frac{1}{\frac{1}{t_1(z-x_1)} -s_2- }\frac{1}{\frac{1}{t_2(z-x_2)}-s_3-}\dots \frac{1}{\frac{1}{t_n(z-x_n)} - h(z)}
\]
where $h$ is an arbitrary function in the Pick class that is meromorphic at $x_1,\dots,x_n$ and satisfies $h(x_j) \neq y_j, \, j=1,\dots,n$.
\end{corollary}
It is of course simple in principle to modify this expression for the case that some $s_k$ are infinite.

We note a relation between the quantities $s_j,t_j$ etcetera that appear in the parametrization theorems.  If one iterates the determinantal formula in Note \ref{detMtilde} one obtains
\[
t_1t_2\dots t_n\prod_{1\leq k<j\leq n}   '(w^k_j-s_k)^2 = \det M
\]
where the prime on the product sign indicates that any factor that is infinite or zero is to be omitted.

\section{Boundary interpolation in the Schur class} \label{schur}
A version of Sarason's Theorem also holds for boundary interpolation by functions in the Schur class.  Its derivation from Theorem \ref{sarason}, with the aid of Cayley transforms, is straightforward and is outlined by Sarason \cite{S}, so we shall merely state the result and go on to describe the parametrization of all solutions, which is easily derived from Theorem \ref{parametrize}.  We denote by $\T$ the unit circle in the complex plane.

\noindent {\bf Problem $\partial NP\Schur$:}   {\em Given distinct points $\xi_1,\dots,\xi_n\in\T$ and given $\eta_1\dots,\eta_n\in\T$ and $\rho_1,\dots,\rho_n\geq 0$, determine whether there is a function $\ph\in\Schur$, analytic at $\xi_1,\dots,\xi_n$, satisfying
\beq\label{mainschur}
\ph(\xi_j) = \eta_j, \qquad (A\ph)(\xi_j) = \rho_j, \qquad j=1,\dots,n,
\eeq
and if there is, describe the class of all such functions $\ph$.

The {\em Pick matrix} of the problem is the matrix $M=[m_{ij}]_{i,j=1}^n$ where
\[
m_{ij}= \left\{ \begin{array}{lll}  \ds \frac{1-\bar\eta_i\eta_j}{1-\bar\xi_i\xi_j} & \mbox{  if  } & i\neq  j, \\ \rho_j & \mbox{  if  } & i=j. \end{array} \right.
\] }

Here $A\ph(\xi)$ denotes (for a function $\ph$ differentiable and non-zero at $\xi=\e^{\ii t}\in\T$) the rate of change of the argument of $\ph(\e^{\ii t})$ with respect to $t$, so that for $\xi\in\T$
\[
(A\ph)(\xi) = \re \frac{\xi\ph'(\xi)}{\ph(\xi)}.
\]
In particular, if $\ph\in\Schur$ and $|\ph(\xi)|=1$ at some $\xi\in\T$ then $|\ph(\e^{\ii t})|^2$ has a local maximum, hence a critical point at $\xi$, from which it follows that $\xi\ph'(\xi)\overline{\ph(\xi)} \in \R$.  Hence, for such $\ph$ and  $\xi$ we have
\beq \label{angderiv}
(A\ph)(\xi) = \frac{\xi\ph'(\xi)}{\ph(\xi)}.
\eeq

\begin{theorem} \label{saraschur}
Problem $\partial NP\Schur$  is solvable if and only if the corresponding Pick matrix is positive definite or minimally positive.  

The solution set of the problem, if non-empty, contains a finite Blaschke product of degree no greater than the rank of the Pick matrix.

The problem is determinate if and only if the Pick matrix is minimally positive.
\end{theorem}

For $\tau\in\T$ we denote by $C_\tau$ the Cayley transform
\[
C_\tau (\la) = \ii \frac{\tau+\la}{\tau-\la}.
\]
$C_\tau$ maps $\D$ to $\Pi$ and $\T$ to $\R \cup\{\infty\}$.
The following is a simple calculation.
\begin{lemma} \label{schurpick}
Let $\ph\in\Schur$, let $\tau,\si, \xi,\eta \in\T$ and suppose that
\[
\bar\tau\xi = \e^{2\ii \psi} \neq 1, \qquad \bar\si\eta = \e^{2\ii\theta} \neq 1.
\]
Suppose that $\ph$ is analytic at $\xi$ and satisfies $\ph(\xi)=\eta$.  If $f=C_\si\circ\ph\circ C_\tau^{-1}$ then $f\in\Pick$, $f$ is analytic at $C_\tau(\xi)$ and
\[
A\ph(\xi) = \frac{\sin^2 \theta}{\sin^2 \psi} f'\circ C_\tau(\xi).
\]
\end{lemma}
\begin{theorem} \label{paraschur}
Suppose that the Pick matrix of Problem $\partial NP\Schur$ is positive definite.  Let
 $\tau, \si \in\T$ be distinct from $\xi_1,\dots,\xi_n$ and $\eta_1,\dots,\eta_n$ respectively.
There exist
\[
s_1,\dots,s_n\in\R\cup\{\infty\}, \qquad t_1,\dots,t_n > 0,  \qquad  y_1,\dots, y_n \in\R\cup\{\infty\}
\]
such that the general solution of Problem $\partial NP\Schur$ is $\ph=C_\si^{-1}\circ f_1\circ C_\tau$ where the functions $f_{n+1},\dots,f_1\in\Pick$ are given recursively by:
\begin{itemize}
\item[\rm(1)] $f_{n+1}$ is any function in $\Pick$ that is meromorphic at $C_\tau(\xi_1),\dots, C_\tau(\xi_n)$ and satisfies $f_{n+1}\circ C_\tau(\xi_j) \neq y_j, \, j=1,\dots, n$;
\item[\rm (2)]  $f_k$ is the augmentation of $f_{k+1}$ at $C_\tau (\xi_k)$ by $s_k, t_k$ for $k=n,\dots,1$.
\end{itemize}
\end{theorem}
\begin{proof}
With the aid of Lemma \ref{schurpick} one sees that a function $\ph\in\Schur$ is a solution of Problem $\partial NP\Schur$ if and only if
$C_\si \circ \ph \circ C_\tau^{-1}$ is a solution of Problem $\partial NP\Pick$ with data
\[
C_\tau(\xi_1),\dots,C_\tau(\xi_n) \in \R, \quad C_\si(\eta_1),\dots, C_\si(\eta_n)\in \R,  \quad   v_1,\dots,v_n       > 0
\]
where, if $\bar\tau\xi_j= \e^{2\ii \psi_j}$ and $\bar\si\eta_j = \e^{2\ii\theta_j}$,
\beq\label{defv}
v_j = \frac{\sin^2 \psi_j}{\sin^2\theta_j} \rho_j, \qquad j=1,\dots,n.
\eeq
By Theorem \ref{parametrize} there exist $s_j \in\R\cup\{\infty\}, t_j > 0, \, y_j \in\R\cup\{\infty\}$ for $j=1,\dots, n$ such that the recursion (1)-(2) parametrizes all solutions of this problem $\partial NP\Pick$.  We then obtain the general solution $\ph$ of Problem $\partial NP\Schur$ by taking $C_\si \circ \ph \circ C_\tau^{-1}=f_1$. \qed
\end{proof}

Note that we can calculate the $s_j, t_j$ and $y_j$ explicitly by substituting $x_j=C_\tau(\xi_j), \, w_j=C_\si(\eta_j)$ and $ v_j$  as given by equation (\ref{defv}) in the formulae in Lemma \ref{inductf}.  We can also interpret the result as a linear fractional or continued fraction parametrization of all solutions, as in Corollaries \ref{linfrac} and \ref{continued}.

\begin{note} \label{reduceschur}
  Reduction and augmentation in the Schur class
\end{note}
It might seem more natural to parametrize the solution set of Problem $\partial NP\Schur$ by elements of $\Schur$ satisfying suitable inequations, rather than elements of $\Pick$.  This is indeed possible, but reduction and augmentation, when transferred to the Schur class by  Cayley transformation, are very cumbersome.  We therefore focus in this paper on the Pick class, even though our original goal was to understand problem $\partial NP\Schur$ and despite the modest complications required to deal with possible infinite target values
(which do not arise in the case of $\Schur$).  Similar considerations apply to boundary interpolation problems in which values of higher derivatives of the function are prescribed.  We show in a subsequent paper \cite{ALY} how Nevanlinna reduction affords an elementary solution of such problems as well.

Jim Agler, Department of Mathematics, University of California at San Diego, La Jolla, CA 92103, U.S.A\\

N. J. Young,  School of Mathematics, Leeds University LS2 9JT, U.K.

\end{document}